\documentclass[12pt,a4paper]{amsart}
\usepackage[latin1]{inputenc}
\usepackage[T1]{fontenc}
\usepackage{geometry}
\usepackage{ifthen}
\usepackage[all]{xy}
\usepackage{enumerate}
\usepackage{amsmath}
\usepackage{amsthm}
\usepackage{amsfonts}
\usepackage{amssymb}
\usepackage{fancyhdr}
\usepackage{ae}
\title{Control theory and the riemann hypothesis:\\[0.1cm]
a roadmap}
\author{Markku Nihtil\"a}
\date{\today}

\begin{document}

\maketitle
\vspace*{-0.5cm}
\begin{center}
{\it University of Kuopio\footnote{{\sc University of Eastern Finland, Department of Physics and Mathematics} (in 2010)}, Department of Mathematics and Statistics\\
POB 1627, FI-70211, Kuopio, Finland}\\[0.2cm]
\texttt{markku.nihtila@uku.fi}

\end{center}
\section{Introduction}
All this seems to be started by Georg Friedrich Bernhard Riemann's original paper whose English translation was included at least in \cite{Edwards74}, without, however forgetting the enormous influence of Leonhard Euler and his studies on the function now known as the Riemann zeta-function.
 
After this, a huge number of papers is available in the internet on the Riemann hypothesis and the zeta-function. Furthermore, several books,  classical, like \cite{Titch86} \& \cite{Edwards74}, and new ones, like \cite{ivic85} \& \cite{borwal09}, have been published on the Riemann zeta-function. The official problem statement on the Riemann hypothesis, which claims that real parts of  the complex zeros of the Riemann zeta-function all are $\frac{1}{2}$, is described in \cite{bombieri01} \& \cite{sarnak02}. A disputed proposal for the proof has been presented \cite{brangesxx} among many erroneous ones.

But the only control-theoretic paper, which we have found, where the Riemann hypothesis and its relation to stability of a given dynamic control system is considered, is \cite{popov85}. On the other hand several authors have considered dynamic systems from the spectral viewpoint relating them to the location of non-trivial, {\it i.e.} complex, zeros of the Riemann zeta-function,  see {\it e.g.} \cite{berrykeating99} \& \cite{connes89}.

There are several statements, which have been proved equivalent to the Riemann hypothesis, see a comprehensive list in \cite{borwal09}. Among them we mention the condition $\Lambda \leq 0$ of the de Bruijin-Newman constant improved by Odlyzko $\Lambda$ \cite{odlyzko00}, Lagarias' statement including harmonic sums \cite{lagarias02}, and Li's condition on the positivity of a certain $\lambda_n$-sequence \cite{borwal09}. The condition, which has a direct connection to our studies concerns the Chebyschev function $\psi(x) =\sum_{p^n < x} \ln p$, where the $p'$s are prime numbers.
It has been proved by von Mangoldt that $\psi$ has the representation,  
see \cite{Edwards74},
\begin{equation}
\psi(x) = x - \sum_{\rho}\dfrac{x^{\rho}}{\rho} - \dfrac{1}{2} \ln (1- x^{-2}) - \ln 2 \pi,
\end{equation}
where the non-trivial complex zeros $\rho$ of the Riemann zeta-function $\zeta(s)$ are counted pairwise in the sum for increasing absolute values of the imaginary parts.
The Riemann hypothesis, represented with the variable $t$ instead of $x=\textrm{e}^t$, is 
equivalent to the statement
\begin{equation}
\psi(\textrm{e}^t) - \textrm{e}^t = O(t^{2} \textrm{e}^{t/2}),
\end{equation}
see  \cite{davenport80}
(The big-O notation is described in \cite{bigO}).
This condition can be considered as a second order asymptote, because the prime number theorem is equivalent to the condition $\psi(x) \sim x$ (asymptotic ratio is equal to 1).

We present here a way to look at the Riemann hypothesis from the viewpoint of Mathematical Control Theory. Firstly, the ratio $G(s)=\frac{1}{(s-1)\zeta(s)}$ is considered as the transfer function (even if transcendental) of a linear dynamic control system. It is developed as a geometric-like series converging around the positive real axis. The individual terms of the series are well-defined transfer functions. These transfer function terms as functions of the complex variable $s$ are inverted into the time domain to obtain so-called impulse responses of the corresponding linear systems. Then these time functions are summed up to obtain the impulse response denoted by $g(t)$ of the system $G(s)$. Our conjecture then concerns the growth bound of the function $g$. The conjecture is: $g(t) = O(t^{k} \textrm{e}^{t/2})$ for some $k \in \{0, 1, 2\}$.
If the conjecture is true then the function $G$ is analytic in 
$\mathbb{C}^+_{1/2 + \varepsilon} = \{s \in \mathbb{C}\, | \, \Re(s) > 1/2+\varepsilon \, \}$ for every $\varepsilon > 0$. This conjecture, if proven, would lead to the truth of the Riemann hypothesis.
\subsection{Simple example 1}
When the transfer function of a linear system is rational, like
$G(s) =\frac{1}{s-a}$, growth of the impulse response $g(t)= \textrm{e}^{ a t}$ and the location of the pole(s) $s_1 =  a$ are directly related to each other. The technique described above for obtaining the impulse response as a series expansion works transparently in this simple example:
\begin{equation} \label{transfer}
G(s) = \dfrac{1}{s-a} = \dfrac{1}{s}\, \dfrac{1}{1- \dfrac{a}{s}} = \sum_{k=0}^{\infty} \dfrac{1}{s}\left(\dfrac{a}{s}\right)^k = \sum_{k=0}^{\infty} \dfrac{a^k}{s^{k+1}}
\end{equation}
\begin{equation} \label{impulse}
g(t) = \mathcal{L}^{-1}\{G(s)\}(t) = \sum_{k=0}^{\infty} \dfrac{1}{k!} (at)^k = \textrm{e}^{at}.
\end{equation}
It is then concluded from the growth condition of Eq. (\ref{impulse}) that the transfer function $G(s)=\frac{1}{s-a}$ is analytic in $\mathbb{C}^+_{a}$.

\section{Main theorems and the growth conjecture}
We put the logic of our reasoning into the form of theorems and lemmas even if some of them are almost immediately evident, if not trivial. The only gap which has not been proved is our growth conjecture, the proof of which may, however, be constructible via modern symbolic program packages.
\newtheorem*{th01}{Theorem 1}
\begin{th01}
The scalar transfer function
\begin{equation}
G(s)=\dfrac{1}{(s-1) \zeta(s)}
\end{equation}
of a linear dynamic system, 
where $\zeta$ is the Riemann zeta-function, representable  e.g. for $\Re(s) > 1$ as the infinite series
\begin{equation}
\zeta(s)=\sum_{m=1}^{\infty} \dfrac{1}{m^s},
\end{equation}
can be represented as a geometric-like series 
\begin{equation}
\begin{split}
G(s)&=\sum_{k=0}^{\infty}\,G_k(s), \\
G_k(s) &= \dfrac{1}{s}\, \Big[1-(s-1)\, \frac{1}{s}\, \zeta (s)\Big]^k, \\
\end{split}
\end{equation}
which converges  at least in a connected set $A \subset \mathbb{C} $ around the positive real axis.
\end{th01}
\begin{proof}
The nominator and denominator of the transfer function can be scaled by such a transfer function, say $R(s)$, which does not change stability properties of $G$ on the strip 
$0 <\Re(s) < 1$. Because we are interested in the location of the poles, or in other words, the location of the zeros of the zeta-function, it is sufficient to use as the scaling function $R(s)=\frac{1}{s}$. Consequently,
\begin{equation} \label{sum01}
\begin{split}
G(s)&=\dfrac{1}{(s-1)\, \zeta (s)} = \dfrac{1}{s}\, \dfrac{1}{(s-1)\,\frac{1}{s}\, \zeta (s)}
= \dfrac{1}{s}\,\dfrac{1}{1-\left[1-(s-1)\, \frac{1}{s}\, \zeta (s)\right]} \\
&=  \dfrac{1}{s}\, \sum_{k=0}^{\infty}\, \Big[1-(s-1)\, \frac{1}{s}\, \zeta (s)\Big]^k \\
&=   \sum_{k=0}^{\infty}\,\dfrac{1}{s}\, \Big[1-(s-1)\, \frac{1}{s}\, \zeta (s)\Big]^k = \sum_{k=0}^{\infty}\,G_k(s).
\end{split}
\end{equation}
The infinite sum converges in the set
$$A = \left\{ s \in \mathbb{C}\,  \Big| \, \Big|1- \dfrac{s-1}{s}\, \zeta(s)\Big| < 1 \, \right\}.$$
The set $A$  is not empty, because 
$\lim_{s \to \infty} Q(s) = 0$ for $s \in \mathbb{R}$ and $Q(s) = 1 -\tfrac{s-1}{s}\, \zeta(s)$ is continuous for $s \in \mathbb{C} \smallsetminus \{0\}$ when suitable forms of analytic continuation of the zeta-function are used.
Then the convergence region can be expanded at least a bit around the positive real axis giving the set $A$.
\end{proof}
\newtheorem*{lm01}{Lemma 1}
\begin{lm01}
The individual terms $G_k(s)$ can be represented as finite sums via the binomial formula
\begin{equation}
G_k(s)= \dfrac{1}{s}\Big[1-(s-1)\, \frac{1}{s}\, \zeta (s)\Big]^k 
= \sum_{n=0}^{k} (-1)^{n}\, \binom{k}{n} \dfrac{(s-1)^{n}}{s^{n+1}}\, \zeta(s)^{n}.
\end{equation}
\end{lm01}
\newtheorem*{th1}{Lemma 2}
\begin{th1}
The inverse Laplace transform
$f_n(t)=\mathcal L ^{-1}\left\{F_n(s)\right\}(t)$ of the function
$$F_n(s) = \dfrac{(s-1)^n}{s^{n+1}}\, \,  \zeta (s)^n,$$

is given, for $0 < t < \ln(N+1)$, by
\begin{equation} \label{formula1}
\begin{split}
f_n(t)& = \mathcal L ^{-1}\left\{\dfrac{(s-1)^n}{s^{n+1}}\, \zeta (s)^n\right\}(t) \\
&= \sum_{m=1}^{N} d_n(m)h(t- \ln m)L_n(t-\ln m),
\end{split}
\end{equation}
where the unit step function is defined by
\begin{equation*}
h(t) =
\begin{cases}
0, & t < 0 \\
\tfrac{1}{2}, & t = 0 \\
1, & t > 0
\end{cases},
\end{equation*}
and $L_n$ is the n'th Laguerre polynomial
$$L_n(t) = \sum_{\nu=0}^{n} \binom{n}{\nu} \dfrac{(-t)^{\nu}}{\nu !} = \dfrac{(-t)^{n}}{n !}+ \cdots - n t+1.$$
The function  $d_n(m)$ is the Piltz divisor function.
\end{th1}
The Piltz divisor function $d_n(m)$, see \cite{Titch86} p. 313,  gives the integer telling how many different ways the integer $m$ can be represented as a product of exactly $n$ integers, {\it i.e.}
$$d_n(m) = \sum_{k_1k_2\cdots k_n=m} \, 1.$$
Some error bounds have been calculated for the summatory divisor function, where the sums of 
$d_n(m)$ were taken over finite number of $m$'s \cite{Olivier}.
The divisor function appears in the powers of zeta-function
$$\zeta(s)^n = \sum_{m=1}^{\infty} \dfrac{d_n(m)}{m^s}.$$

\vspace*{0.2cm}
\noindent{\it Proof of Lemma 2.}
The terms $\dfrac{1}{m^s} = \textrm{e}^{- \ln m \, s}$ are Laplace transforms of delta-distributions delayed by \mbox{$\tau_m =\ln m$}, {\it i.e.}
$$\mathcal L ^{-1}\left\{\zeta(s)^n\right\}(t)= \sum_{m=1}^{N}d_n(m) \delta (t-\ln m), \ \ 0 < t < \ln(N+1),$$
and Laguerre polynomials are obtained as
$$L_n(t)=\mathcal L ^{-1}\left\{\dfrac{(s-1)^n}{s^{n+1}}\right\}(t).$$
Because the function $f_n(t)=\mathcal L ^{-1}\left\{F_n(s)\right\}(t)$ can be represented as a sum of convolution integrals of delayed delta-functions and the Laguerre polynomial $L_n$ the result (\ref{formula1}) is immediate.\hspace*{\stretch{1}}$\square$

\vspace*{0.2cm}
The preceeding results are collected finally in the form of a theorem the proof of which is immediate.
\newtheorem*{th02}{Theorem 2}
\begin{th02}
The inverse transforms of the individual terms $g_k(t)= \mathcal{L}^{-1}\{G_k(s)\}(t)$ of the expansion in (\ref{sum01}) are
\begin{equation} \label{geekoo}
g_k(t) = \sum_{n=0}^{k}(-1)^n \binom{k}{n} f_n(t)
 = \sum_{n=0}^{k}(-1)^n \binom{k}{n} \sum_{m=1}^{N} d_n(m)h(t- \ln m)L_n(t-\ln m),
\end{equation}
and the inverse Laplace transform (desired impulse response) $g(t) = \mathcal{L}^{-1}\{G(s)\}(t)$ 
for $0 < t < \ln (N+1)$, assuming that the series converges, is
\begin{equation} \label{response}
g(t) = \sum_{k=0}^{\infty} \sum_{n=0}^{k}(-1)^n \binom{k}{n} 
\sum_{m=1}^{N} d_n(m)h(t- \ln m)L_n(t-\ln m).
\end{equation}
\end{th02}
For $\ln N < t$ in (\ref{geekoo}) \& (\ref{response}) the step functions $h(t- \ln m)$ can be dropped out, because they all are equal to 1 \cite{Olivier02}. In order to keep the sum (\ref{response}) in a concise form we have counted the terms from $n=0$. Actually, then the coefficient  $d_0(m)$ is not defined in number theory. This discrepancy can be avoided by defining artificially this coefficient. At present, it seems that the correct definition is $d_0(m) =1$. 
\newtheorem*{th03}{The growth conjecture}
\begin{th03}
The impulse response $g(t)$, Eq.(\ref{response}), of the system whose transfer function is $G(s)= \frac{1}{(s-1) \zeta(s)}$ has a limited growth in the sense that for some $k \in \{0, 1, 2\}$
$$g(t)= O(t^k \textrm{e}^{t/2}).$$
\end{th03}
\newtheorem*{th13}{Theorem 3}
\begin{th13}
If the growth conjecture is true, then also the Riemann hypothesis is true, i.e. all the complex zeros of  the Riemann zeta-function $\zeta(s)$ are of the form $s_k = \frac{1}{2}\pm \mathfrak{i}\, \gamma_k$ where $ \gamma_k \in \mathbb{R}$.
\end{th13}
\newtheorem*{lm03}{Lemma 3\footnote{Results of this Lemma were kindly pointed out by Hans Zwart, University of Twente, \\ \hspace*{0.6cm}The Netherlands}}
\begin{lm03}
(c.f. Property A.6.2, \cite{swart95})  Laplace transformable functions 
\mbox{$h:[0, \infty) \to \mathbb{R}$} have the property: If e$^{\,- \beta t} h(t) \in {\bf L}_1([0, \infty); \mathbb{R})$  for some real $\beta$, then 
$H(s)=\mathcal{L}\{h(t)\}(s)$ is holomorfic (analytic) and bounded on $\mathbb{C}^+_{\beta}=\{s \in \mathbb{C} \, | \, \Re(s) > \beta \} $.
\end{lm03}
\noindent{\it Proof of Theorem 3.} The growth conjecture gives  
the bound $|g(t)| \leq M t^k \textrm{e}^{t/2}$ for some $M \in \mathbb{R}$. The impulse response $g(t)$ is continuous except a countable number of points. Then $g(t)$ is locally Lebesgue integrable and  
$\textrm{e}^{-(1/2 + \varepsilon)t} g(t) \in {\bf L}_1([0, \infty); \mathbb{R})$ for all $\varepsilon > 0.$
Then, based on Lemma 3, $G(s)$ is analytic on $\mathbb{C}^+_{1/2+\varepsilon}$ for all $\varepsilon > 0$. Due to symmetry of the complex poles of $G(s)$ (or zeros of the zeta-function) with respect to the line $\Re(s) =\frac{1}{2}$, poles of $G(s)$ cannot either be on 
$\mathbb{C}^-_{1/2 -\varepsilon}=\{s \in \mathbb{C} \, | \, \Re(s) < \frac{1}{2} -\varepsilon \} $ for any $\varepsilon > 0.$ Consequently, the complex poles of $G(s)$, and correspondingly the complex zeros of the Riemann zeta-function $\zeta(s)$ must be on the line 
$\Re(s) =\frac{1}{2}$ of the complex plane $\mathbb{C}$.\hspace*{\stretch{1}}$\square$
\section{Conluding remarks}
\noindent {\bf Remark 1.} It has to be noted that the implication of Lemma 3 cannot be reversed for general transcendental transfer functions. Also only partial results can be obtained, see 
\cite{macc00}. In the case of rational transfer functions the location of poles and growth conditions only are equivalent in a certain sense.

\vspace*{0.2cm}
\noindent {\bf Remark 2.} The formula (\ref{response}) gives the third, and new, pair of a function and its Laplace transform where the transform has as its poles the zeros of the Riemann zeta-function, {\it i.e.}
$$\mathcal{L}^{-1}\left\{\dfrac{1}{(s-1)\, \zeta(s)}\right\}(t) = g(t).$$
{ }\\
The two other pairs, well-known by number theorists,  are
$$\mathcal{L}^{-1}\left\{-\, \dfrac{\zeta'(s)}{s \, \zeta(s)}\right\}(t) = \sum_{p^n < \textrm{e}^t} \ln p,$$
$$\mathcal{L}^{-1}\left\{-\, \dfrac{\zeta'(s)}{s \, \zeta(s)}\right\}(t) 
= \textrm{e}^t - \sum_{\rho} \dfrac{\textrm{e}^{\rho t}}{\rho}\, - \,  \dfrac{1}{2} \, \ln \big( 1- \textrm{e}^{-2t} \big) -\ln 2 \pi. $$
The last pair is a result of the inverse Mellin transform, see \cite{Flajolet95} \& \cite{Titch86}, for $x= \textrm{e}^t$.

\vspace*{0.2cm}
\noindent{\bf Remark 3.} It is certainly temptating to try to calculate the series (\ref{response}) by using some symbolic programs to see if the growth conjecture looks (maybe is) true.

\vspace*{0.2cm}
\noindent{\bf Remark 4.}  It is known that $|L_n(t)| \leq \textrm{e}^{t/2}$ for all $t> 0$. 
Then a possible road to the truth of the growth conjecture is to find (maybe analytically) for $L_n$ an expression like
$L_n(t) = \pm \textrm{e}^{t/2}$+asymptotically small error term.


\begin{thebibliography}{10}

\bibitem{berrykeating99}
Berry, M.V. and Keating, J.P., The Riemann zeros and eigenvalue asymptotics, {\it SIAM Review}, 
Vol. 41, No. 2, 1999, pp. 236--266.

\bibitem{bigO}
Big O notation: available in 
\texttt{http://en.wikipedia.org/wiki/Big\underline{ }O\underline{ }notation}
\bibitem{bombieri01}
Bombieri, E., Problems of the millennium: The Riemann hypothesis, Preprint, 11 pp.
(available in: \texttt{http://www.claymath.org/millennium})


\bibitem{Olivier}
Bordell\`es, O., Explicit upper bounds for the average order of $d_n(m)$ and application to class number, {\it Journal of Inequalities in Pure and Applied Mathematics}, Vol. 3, No. 3, Article 38, 2002, 15 pp.

\bibitem{Olivier02}
Bordell\`es, O., Personal communication, March 2009.

\bibitem{borwal00}
Borwein, J.M., Bradley, D.M., and Crandall, R.E., Computational strategies for the Riemann
 zeta function, {\it Journal of Computational and Applied Mathematics}, Vol. 121, 2000, pp. 247--296.

\bibitem{borwal09}
Borwein, P., Choi, S., Rooney, B., and Weirathmueller, A. (Eds.), {\it The Riemann Hypothesis, A resource for the Afficionado and Virtuoso Alike}, Springer, New York, N.Y., U.S.A., 2008.

\bibitem{brangesxx}
de Branges, L., A proof of the Riemann hypothesis, Preprint, July 19  2006, 41 pp.

\bibitem{connes89}
Connes, A., Trace formula in noncommutative geometry and the zeros of the Riemann zeta function, 
{\it Sel. Math. New ser.}, Vol. 5, No. 1, 1999, pp. 29--106. \\
(also available in: \texttt{http://www.alainconnes.org/en/downloads.php})

\bibitem{swart95}
Curtain, R.F. and Zwart, H.J., {\it An Introduction to Infinite-Dimensional Linear Systems Theory}, New York, Springer, 1995.

\bibitem{davenport80}
Davenport, H., {\it Multiplicative Number Theory}, 2nd ed., New York, U.S.A., Springer, 1980, p. 104.

\bibitem{Edwards74}
Edwards, H.M., {\it Riemann's Zeta Function}, Mineola, N.Y., U.S.A., Dover Publications Inc., 2001.

\bibitem{Flajolet95}
Flajolet, P., Gourdon, X., and Dumas, P., Mellin transforms and asymptotics: Harmonic sums, {\it Theoretical Computer Science}, Vol. 144, 1995, pp. 3--58
(available in: \texttt{http://algo.inria.fr/flajolet/Publications/FlGoDu95.pdf})

\bibitem{ivic85}
Ivi\'c, A., {\it The Riemann Zeta-Function, Theory and Applications}, 
Mineola, N.Y., U.S.A., Dover Publications Inc., 2003.

\bibitem{lagarias02}
Lagarias, J.C., An elementary problem equivalent to the Riemann hypothesis, 
{\it The American Mathematical Monthly}, Vol. 109, No. 6, 2002, pp. 534--543.

\bibitem{macc00}
MacCluer, C.R., An unstable plant with no poles, {\it IEEE Trans. Automatic Control}, Vol. 45, No. 8, Aug. 2000, pp. 1575--1576.

\bibitem{odlyzko00}
Odlyzko, A.M., An improved bound for the Bruijin-Newman constant, {\it Numerical Algorithms}, Vol. 25, 2000, pp. 293--303.

\bibitem{popov85}
Popov, V.M., On stability properties which are equivalent to Riemann hypothesis, {\it Libertas Mathematica}, Vol. 5, 1985, pp. 55--61.

\bibitem{sarnak02}
Sarnak, P., Problems of the millennium: The Riemann hypothesis (2004), Preprint, 9 pp. 
(available in: \texttt{http://www.claymath.org/millennium})

\bibitem{Taylor07}
Taylor, J. L., {\it Complex Variables}, Version 2.1., 2007. (available in: \\ \texttt{http://www.math.utah.edu/$\sim$taylor/9\underline{ }GammaZeta.pdf})

\bibitem{Titch86}
Titchmarsh, E.C., {\it The Theory of the Riemann Zeta-function}, 2nd ed. (Revised with comments by D.R.~Heath-Brown), Oxford, U.K., Oxford University Press, 1986.

\end{thebibliography}
\end{document}